\documentstyle{article}

\newtheorem{defn}{Definition}[section]
\newtheorem{lem}[defn]{Lemma}
\newtheorem{thm}[defn]{Theorem}
\newtheorem{prop}[defn]{Proposition}
\newtheorem{cor}[defn]{Corollary}
\newtheorem{rem}[defn]{Remark}
\newtheorem{ex}[defn]{Example}
\newtheorem{assu}[defn]{Assumption}
\newtheorem{constr}[defn]{Construction}
\makeatletter
\def\theequation{\thedefn.\@arabic\c@equation}
\makeatother

\newcommand{\A}{{\cal A}}

\newcommand{\D}{{\cal D}}
\newcommand{\X}{{\cal X}}
\newcommand{\Y}{{\cal Y}}
\newcommand{\HH}{{\cal H}}
\newcommand{\M}{{\cal M}}
\renewcommand{\O}{{\cal O}}
\renewcommand{\L}{{\cal L}}

\newcommand{\W}{{\cal W}}
\newcommand{\PP}{{\cal P}}

\newcommand{\VV}{{\cal V}}

\newcommand\ZZ{{\cal Z}}
\newcommand{\C}{{\bf C}}
\renewcommand{\P}{{\bf P}}
\newcommand{\Q}{{\bf Q}}
\newcommand{\Z}{{\bf Z}}

\newcommand{\De}{\Delta}
\renewcommand{\phi}{\varphi}
\newcommand{\eps}{\varepsilon}
\newcommand{\epsi}{\eps^i_{\chi,\chi'}}
\newcommand{\into}{\hookrightarrow}
\newcommand{\<}{\langle}
\renewcommand{\>}{\rangle}

\newcommand{\mod}{\hbox{\rm mod}\;}
\newcommand{\longdownarrow}{\downarrow}
\newcommand{\isom}{=}
\newcommand{\NS}{\hbox{\sl NS}}

\newcommand{\ord}[1]{{\rm o}({#1})}
\def\re^#1_#2{a^{#1}_{#2}}
\def\reb^#1_#2{(n_{#2}a^{#1}_{#2})/m_{#1}}

\newcommand{\autg}{\Phi}

\newcommand{\ichi}{{i,\chi}}
\newcommand{\defor}[1]{\hbox{\rm Def}_{#1}}
\newcommand{\defg}[1]{\hbox{\rm Def}^G_{#1}}
\newcommand{\defnat}{\hbox{\rm Dnat}_X}
\newcommand{\defgal}{\hbox{\rm Dgal}_X}

\newcommand{\hilbd}{Hilb^{\rm div}}
\newcommand{\birto}{\hbox{\rm-}\,\hbox{\rm-}\,\hbox{\rm-}\!\!\!\!\!>}
\newcommand{\Ger}{\hbox{\sl Germs}}
\newcommand{\Sets}{\hbox{\sl Sets}}

\newcommand{\Pointed}{\hbox{\sl Ansp}_0}
\newcommand{\Art}{\hbox{\sl Art}^*}
\newcommand{\Pf}{\noindent{\sc Proof.~}}
\newcommand{\qed}{\ $\Box$\par\smallskip}

\newcommand\Bigoplus{\textstyle\bigoplus\limits}
\newcommand\Bigotimes{\textstyle\bigotimes\limits}

\newcommand{\al}{\alpha}
\newcommand{\be}{\beta}
\newcommand{\albe}{{\al\be}}
\newcommand{\beal}{{\be\al}}
\newcommand{\si}{\sigma}
\newcommand{\alga}{{\al\ga}}
\newcommand{\bega}{{\be\ga}}
\newcommand{\ga}{\gamma}
\newcommand{\la}{\lambda}

\title{Automorphisms and moduli spaces\\
of varieties with ample canonical class\\
via deformations of abelian covers}
\author{Barbara \hbox{\rm Fantechi${}^*$} \quad --- \quad Rita Pardini
\thanks{Both authors are members of GNSAGA of CNR.}}
\date{}

\begin{document}
\maketitle
\begin{abstract}By a recent result of Viehweg, projective manifolds with
ample canonical
class have a coarse moduli space, which is a union of quasiprojective
varieties. In this paper, we prove that there are manifolds with ample
canonical class that lie on arbitrarily many irreducible components of the
moduli; moreover, for any finite abelian group $G$ there exist infinitely
many components $M$ of the moduli of varieties with ample canonical class
such that the generic automorphism group $G_M$ is equal to $G$.
\par
In order to construct the examples, we use abelian covers. Let $Y$ be a
smooth complex projective variety of dimension
$\ge 2$. A Galois cover $f:X\to Y$ whose Galois group is finite and abelian
is called an abelian cover of $Y$; by \cite{Pa1}, it is determined by its
building data, i.e.~by the branch divisors and by some line bundles on $Y$,
satisfying appropriate compatibility conditions. Natural
deformations of an abelian cover are also introduced in \cite{Pa1}.
\par In this paper we prove two results about abelian covers: first, that if
the building data are sufficiently ample, then the natural deformations
surject on the Kuranishi family of $X$; second, that if the building data
are sufficiently ample and generic, then $Aut(X)=G$.
\par These results, although in some sense ``expected'', are in fact rather
powerful and enable us to construct the required examples. Finally, note
that it is essential for our applications to be able to deal with general
abelian covers and not only with cyclic ones.
\end{abstract}

\section{Introduction} Coverings of algebraic varieties are a classical theme
in algebraic geometry, since Riemann's description of curves as branched
covers of the projective line. Double covers were used by the Italian school
to construct examples that shed light on the theory of surfaces and to
describe special classes of surfaces, as in the case of Enriques surfaces.

More recently, cyclic coverings have been extensively applied by several
authors to the study of surfaces of general type; it will be enough to recall
the work of Horikawa, Persson and Xiao Gang. Abelian covers have been used
by Hirzebruch to give examples of surfaces of general type on and near the
line $c_1^2=3c_2$; Catanese and Manetti have used bidouble and iterated
double covers, respectively, of $\P^1\times \P^1$ to construct explicitly
connected components of the moduli space of surfaces of general type.

In \cite{Pa1}, the second author has given a complete description of abelian
covers of algebraic varieties in terms of the so-called building data, namely
of certain line bundles and divisors on the base of the covering, satisfying
suitable compatibility relations. Natural deformations of an abelian cover
$f:X\to Y$ are also introduced there and it is shown that they are complete,
if $Y$ is rigid,  regular and of dimension $\ge 2$, and if the building data
are sufficiently ample. (Natural deformations are obtained by modifying the
equations defining $X$ inside the total space of the bundle $f_*\O_X$).

In this paper we study natural deformations of an abelian cover $f:X\to Y$
and prove that they are complete for varieties of dimension at least two if
the branch divisors are sufficiently ample. The result requires no
assumption on
$Y$, and in particular also holds when the cover has obstructed
deformations; this is a key technical step towards the moduli space
constructions described below.

We then turn to the study of the automorphism group of the cover. Since the
automorphism group of a variety of general type is finite, one would expect
that in the case of a Galois cover it coincides with the Galois group, at
least if the cover is generic. Our main theorem \ref{mainthm} shows that this
is indeed the case for an abelian cover, if the branch divisors are generic
and sufficiently ample.

We construct explicitly coarse moduli spaces of abelian covers and
complete families of natural deformations for a fixed base of the cover $Y$;
this is useful if one wants to investigate the birational structure of the
components of the moduli obtained by the methods of this paper.

The main application of the results described so far is the study of moduli
of varieties with ample canonical class. Recently Viehweg proved the
existence of a coarse moduli space for varieties with ample canonical class
of arbitrary dimension, generalizing Gieseker's result for surfaces. Given
an irreducible component
$M$ of the moduli space of varieties with ample canonical class, the
automorphism group $G_M$ of a generic variety in $M$ is well-defined. In
contrast with the case of curves (where this group is trivial for $g\ge 3$),
it was already known in the case of surfaces that there exist infinitely
many components
$M$ of the moduli with nontrivial automorphism group $G_M$; it is easy to
construct examples such that $G_M$ contains an involution, and Catanese
gave examples where $G_M$ contains a subgroup isomorphic to $\Z_2\times
\Z_2$. There are also, of course, easy examples of components $M$ where
$G_M$ is trivial (for instance the hypersurfaces of degree $d\ge 5$ in $\P^3$).

As a first application of theorem \ref{mainthm} we prove that for any finite
abelian group $G$ there are infinitely many irreducible components $M$ of the
moduli of varieties with ample canonical class such that $G_M=G$; notice that
we precisely determine $G_M$ instead of just bounding it from below.

We also prove that there are varieties with ample canonical class lying on
arbitrarily many irreducible components of the moduli. We distinguish these
components by means of their generic automorphism group; there are examples
both in the equidimensional and in the non-equidimensional case. In the
surface case, this answers a question raised by Catanese in \cite{Ca2}.

Let $S$ be a surface of general type; Xiao has given explicit upper bounds
both for the cardinality of $Aut(S)$ and of an abelian subgroup of $Aut(S)$,
in terms of the invariants of $S$ (\cite{Xi1}, \cite{Xi2}). Some upper
bounds are also known for a higher-dimensional variety $X$ with
ample canonical class, although sharp bounds are still
lacking. It seems interesting to ask whether these bounds can be improved by
considering instead of $Aut(X)$ the group
$Aut_{\rm gen}(X)$, namely the intersection in
$Aut(X)$ of the images of the generic automorphism groups $G_M$ of all
irreducible components $M$ of the moduli space containg $X$ (in particular,
if $X$ lies in a unique component $M$, then $Aut_{\rm gen}(X)=G_M$).

As a first step towards the computation of a sharp bound for
$\#Aut_{\hbox{\rm gen}}(S)$, we show that such a bound cannot be ``too
small''; in fact we give a sequence of surfaces $S_n$ of general type, whose
Chern numbers tend to infinity with $n$, and such that
$\#Aut_{\hbox{\rm gen}}(S_n)\ge 2^{-4}K_{S_n}^2$.

The paper goes as follows: in section 2 we collect some results from the
literature and set up the notation. In section 3 we prove that, if the
branch divisors are
sufficiently ample, then infinitesimal natural
deformations are complete. In section 4 we prove (theorem \ref{mainthm}) that
the
automorphism group of an abelian cover coincides with the Galois group if the
building data are sufficiently ample and generic. To do this, we prove some
results on extensions of automorphisms, which we believe should be of
independent interest. The proof of \ref{mainthm} is based on a degeneration
argument and requires an explicit partial desingularization, contained in
section 7. Section 5 contains the construction of a
coarse moduli space for abelian covers of a given variety $Y$ and of a
complete family of natural deformations. Finally, in section 6 we apply
the results of sections 3 and 4 to
the study of moduli spaces of varieties with ample canonical class, as stated
above.

\smallskip
\noindent{\em Acknowledgements}. This work was supported by the italian MURST
60\% funds. The first author would also like to thank the Max-Planck-Institut
f\"ur Mathematik (Bonn) for hospitality and the italian CNR for support.
\section{Notation and conventions} All varieties will be complex, and smooth
and projective unless the contrary is explicitly stated.

For a projective morphism of schemes $Y\to S$, $Hilb_S(Y)$ will be the
relative Hilbert scheme (see \cite{FGA}, expos\'e 221). When $Y$ is smooth
over $S$, $Hilb^{\rm div}_S(Y)$ will be the (open and closed) subscheme of
$Hilb_S(Y)$ parametrizing divisors (see \cite{Fo} for a proof of this).
When $S$ is a point, it will be omitted
from the notation.

For $Y$ a smooth
projective variety, let
$c_1:Pic(Y)\to H^2(Y,\Z)$ be the map associating to a line bundle its first
Chern class; let
$NS(Y)$ be the image in
$H^2(Y,\Z)$ of
$Pic(Y)$, and
$Pic^\xi(Y)$ the inverse image of $\xi\in NS(Y)$. Let $q(Y)=\dim
H^1(Y,\O_Y)$ be the dimension of $Pic^0(Y)$.

Let $\X\to B$ be any flat family, with integral fibres. Then there are open
subschemes $Aut_{\X/B}$ and $Bir_{\X/B}$ of the relative Hil\-bert sche\-me
\hbox{$Hilb_B(\X\times_B\X)$} parametrizing fibrewise the (graphs of)
automorphisms and birational automorphisms of the fibre (\cite{FGA},
\cite{Ha}).

We denote the cardinality of a (finite) set $S$ by $\#S$; for each integer
$m\ge 2$, let $\zeta_m=e^{2\pi i/m}$.

\smallskip
\noindent{\em Notation for abelian covers}. The following notation will be
used freely throughout the paper: we collect it here for the reader's
convenience.

$G$ will be a finite abelian group, $G^*$ its dual; the
order of an element $g$ will be denoted by $\ord{g}$.
Let $I_G$ be the set of all pairs $(H,\psi)$ where $H$ is a cyclic subgroup
of $G$ with at least two elements and $\psi$ is a generator of $H^*$. There
is a bijection between $I_G$ and $G\setminus 0$ given by $(H,\psi)\mapsto g$
where $g\in H$ is such that $\psi(g)=\zeta_{\#H}$. For $\chi\in G^*$,
$i=(H_i,\psi_i)\in I_G$, let $\re^i_\chi$ be the unique integer such that
$0\le \re^i_\chi<m_i$ (where $m_i=\#H_i$) and
$\chi_{|H_i}=\psi_i^{\re^i_\chi}$ (cfr.\ \cite{Pa1}, remark 1.1  on p.~195,
where $\re^i_\chi$ is denoted by $f_{H,\psi}(\chi)$). Let
$\epsi=[(\re^i_\chi+\re^i_{\chi'})/m_i]$, where $[r]$ is the integral part of
a real number $r$; note that $\epsi$ is either $0$ or $1$.

A basis of $G$ will be a sequence of elements of $G$, $(e_1,\dots,e_s)$,
such that $G$ is the direct sum of
the (cyclic) subgroups generated by the $e_j$'s, and such that $\ord{e_j}$
divides $\ord{e_{j+1}}$ for each $j=1,\ldots,s-1$.
Given a basis
$(e_1,\dots,e_s)$ of $G$,  we will call dual basis of
$G^*$ the $s$-tuple
$(\chi_1,\ldots,\chi_s)$, where $\chi_j(e_i)=1$ if $i\ne j$ and
$\chi_i(e_i)=\zeta_{\ord{e_i}}$.
We will write $\re^i_j$ instead of $\re^i_{\chi_j}$, for all
$j=1,\ldots,s$; for $\chi=\chi_1^{\al_1}\cdots\chi_s^{\al_s}$, let
$$q^i_\chi=\left[\sum_{j=1}^s\frac{\al_i\re^i_j}{m_i}\right].$$
Note that, unlike $\re^i_\chi$, $q^i_\chi$ depends on the choice of the basis
and not only on $\chi$ and $i$.

\begin{lem} Let $G$ be as above, and let $I\subset I_G$ be a subset with $k$
elements {\rm(}which we denote by $1,\ldots,k${\rm)} such that  the natural
map $H_1\oplus
\ldots \oplus H_k\to G$ is surjective. Then the $k\times s$ matrix
$(\re^i_j)$ has rank $s$ over $\Q$.
\end{lem}
\Pf
Let $g_i$ be the element corresponding to $(H_i,\psi_i)$ via the bijection
$I_G\leftrightarrow G\setminus 0$ described above. Then, for any
$i=1,\ldots,k$ and for any $j=1,\ldots,s$, one has $\re^i_j/m_i=\la_{ij}/n_j$,
where $n_j=\ord{e_j}$ and $g_i=\sum \la_{ij}e_j$, with $0\le \la_{ij}<n_j$ and
$\la_{ij}\in\Z$ . So the matrix $(\re^i_j)$ has the same rank over $\Q$ as the
matrix $\la_{ij}$. On the other hand $\la_{ij}$ is the matrix associated to
the natural map $H_1\oplus
\ldots \oplus H_k\to G$, which is surjective. Let $p$ be a prime factor of
$n_1$, hence of all of the $n_j$'s. Then the map $\Z_p^k\to \Z_p^s$
represented by the matrix $(\la_{ij})\ \mod p$ is also surjective, hence the
matrix $(\la_{ij})$ has an $s\times s$ minor whose determinant is nonzero
modulo $p$. This implies that the determinant is nonzero, hence the result.
\qed
\smallskip

Let $X$ be any projective variety. A {\sl deformation} of $X$ over a pointed
analytic space
$(T,o)$ will be a flat, proper map $\X\to T$, together with an isomorphism of
the special fibre $\X_o$ with $X$.

Deformations modulo isomorphism are a
contravariant functor $Def_X$ from the category $\Pointed$ of pointed analytic
spaces to the category $\Sets$, where the functoriality is given by pullback.

More generally, given a contravariant functor $F:\Pointed\to \Sets$, we will
use the same letter $F$ to denote the induced functor on the categories $\Ger$
of germs of analytic spaces and $\Art$ of finite length spaces supported in a
point (i.e. $Spec$'s of local Artinian $\C$-algebras). For the properties of
functors on $\Art$, we refer the reader to \cite{schl}.

\smallskip

Let $M$ be an irreducible component of the moduli space of
(projective) manifolds with ample canonical class. As the automorphism group
is semicontinuous (see corollary 4.5), it makes sense to speak of the
automorphism group of a generic manifold in $M$; we will denote it by $G_M$.
Note that for any $X$ such that $[X]\in M$, there is a natural identification
of $G_M$ with a subgroup of $Aut(X)$. If $X$ is a minimal surface of general
type, we denote the intersection in
$Aut(X)$ of $G_M$ for all components $M$ containing $[X]$ by $Aut_{\hbox{\rm
gen}}(X)$; it is the largest subgroup $H$ of $Aut(X)$ such that the action of
$H$ extends to any small deformation of $X$.

\section{Deformations of abelian covers}

In this section we
introduce natural deformations of a smooth abelian cover and prove
that infinitesimal natural deformations
are complete, if  the branch
divisors are sufficiently ample and the dimension is at least two.

We start by recalling from \cite{Pa1} some fundamental results on abelian
covers; the reader will find there a more detailed exposition and proofs of
the following statements.

Let $G$ be a finite abelian group and
let
$I$ be a subset of
$I_G$: we will use freely throughout the paper the notation introduced in
section 2. Let
$Y$ be a smooth projective variety: a $(G,I)$-cover of
$Y$ is a normal variety $X$ and a Galois cover $f:X\to Y$ with
Galois group $G$ and branch divisors $D_i$ (for $i\in I)$ having
$(H_i,\psi_i)$ as inertia group and induced character (see \cite{Pa1} for
details). $X$ is smooth if and only if the $D_i$'s are smooth, their union is
a normal crossing divisor, and, whenever
$D_{i_1},\ldots,D_{i_k}$ have a common point, the natural map
$H_{i_1}\oplus\ldots\oplus H_{i_k}\to G$ is injective. The cover is said to be
{\em totally ramified} if the natural map $\bigoplus_{i\in I}H_i\to G$ is
surjective. Note that each abelian cover can be factored as the composition of
a totally ramified with an unramified cover.

Let $M_i=\O_Y(D_i)$. The vector bundle $f_*\O_X$ on $Y$ splits naturally as sum
of eigensheaves
$L_\chi^{-1}$ for $\chi\in G^*$, and multiplication in the $\O_Y$-algebra
$f_*\O_X$ induces isomorphisms \begin{equation}
\label{bdata}L_\chi\otimes L_{\chi'}=L_{\chi\chi'}\otimes\Bigotimes_{i\in
I} M_i^{\otimes \epsi}\qquad\quad \hbox{for all
$\chi,\chi'\in G^*\setminus 1$}.
\end{equation}

Denote $L_{\chi_j}$ by $L_j$, and let $n_j=\ord{\chi_j}$. The isomorphisms
above induce isomorphisms

\begin{equation} \label{rbdata} L_j^{\otimes n_j} =\Bigotimes_{i\in
I}M_i^{\otimes\reb^i_j}\qquad\quad \hbox{for all $j=1,\ldots,s$}.
\end{equation}
The $(D_i,L_\chi)$ are the {\em building data} of the cover; the $(D_i,L_j)$
are the {\em reduced building data}. The sheaves $L_\chi$ can be recovered from
the reduced building data by setting, for
$\chi=\chi_1^{\al_1}\cdots\chi_s^{\al_s}$,
\begin{equation} \label{chidarbd}
L_\chi=\Bigotimes_{j=1}^s L_j^{\alpha_j}\otimes\Bigotimes_{i\in
I}M_i^{-q^i_\chi}.
\end{equation}
Conversely, for each choice of $(D_i,L_\chi)$ (resp.\ $(D_i,L_j)$) satisfying
equation (\ref{bdata}) (resp.\ (\ref{rbdata})), there exists a unique cover
having these as (reduced) building data. Note that equations (\ref{rbdata})
have
a solution in
$Pic(Y)$ (viewing the line bundles $M_i$'s as parameters and the
$L_j$'s as variables) if and only if their images via $c_1$ have a solution
in $\NS(Y)$.

\begin{assu}\label{totram} In this paper all $(G,I)$-covers will be totally
ramified.
Unless otherwise stated, $f:X\to Y$ will be a
$(G,I)$-cover, with reduced building data $(D_i,L_j)$. We will also assume
that
$X$ and
$Y$ are smooth, of dimension $\ge 2$, and that $X$ has ample canonical class.
\end{assu}
We say that a property holds whenever a line bundle $L$ (or a divisor $D$) is
sufficiently ample if it holds whenever $c_1(L)$ (or $c_1(D)$) belongs to a
(given) suitable translate of the ample cone. It is easy to see that
assumption \ref{totram} implies the following: if all of the $D_i$'s are
sufficiently ample then so is $L_\chi$ for any $\chi\ne 1$.  Moreover, if $V$
is a vector bundle, $V\otimes L$ is ample for any sufficiently ample $L$.

Let $S=\{(\ichi)\in I\times G^*|\chi_{|H_i}\ne\psi_i^{-1}\}$. Given a
$(G,I)$-cover
$X\to Y$ as above, together with sections
$s_\ichi$ of
$H^0(M_i\otimes L_\chi^{-1})$ for all $(\ichi)\in S$, a natural
deformation of
$X$ was defined in
\cite{Pa1}, \S 5.
We now give a functorial (and more general) version of that definition in
order to be able to apply standard techniques from deformation theory.

\begin{defn}\label{natdef}{\rm A {\em natural deformation of the reduced
building data} of
$f:X\to Y$ over $(T,o)\in\Pointed$ is $(\Y,\M_i,\L_j,s_{\ichi},\phi_j)$
where:\begin{enumerate}
\item $i\in I$, $j=1,\ldots,r$, and $(\ichi)\in S$;
\item $\Y\to T$ is a deformation of $Y$ over $T$;
\item $\L_j$ and $\M_i$ are line bundles on $\Y$
such that $\L_j$ restricts to $L_j$ and $\M_i$ to $M_i$ over $o$;
\item $
\phi_j:\L_j^{\otimes n_j}\to \bigotimes \M_i^{\otimes \reb^i_j}$ is an
isomorphism whose restriction to $\Y_o$ coincides with the isomorphism
$L_j^{\otimes n_j}\to \bigotimes M_i^{\otimes \reb^i_j}$ given by
multiplication;
\item $s_{\ichi}$ is a section of $\L_\chi^{-1}\otimes\M_i$, where $
\L_\chi=\Bigotimes_{j=1}^s \L_j^{\alpha_j}\otimes\Bigotimes_{i\in
I}\M_i^{-q^i_\chi}$;
\item $s_{\ichi}$ restricts over $\Y_o$ to $s_\ichi^0$, where
$s_\ichi^0=0$ if
$\chi\ne 1$, and $s^0_{i,1}$ is a section of $M_i$ defining
$D_i$.
\end{enumerate}
We will say that a deformation is {\em Galois} if $s_{i,\chi}=0$ for $\chi\ne
1$.}
\end{defn}

Natural deformations modulo isomorphism define a contravariant functor
$\defnat:\Pointed\to
\Sets$, and Galois deformations are a subfunctor $\defgal$. Note that the
inclusion
$\defgal\hookrightarrow \defnat$ is naturally split.
We now extend formulas in \S 5 of \cite{Pa1} to define a natural
transformation of functors
$\defnat\to \defor X$.
\begin{defn}\label{trasnat}{\rm Let $T$ be a germ of an analytic space, and let
$$(\Y,\L_j,\M_i,\phi_j,s_\ichi)\in\defnat(T).$$ Let $V$ be the
total space of the vector bundles $\bigoplus_{\chi\in G^*}\L_\chi$, and let
$\pi:V\to \Y$ be the natural projection. For a line bundle $\L$ on $\Y$,
denote its pullback to $V$ by $\bar \L$, and analogously for sections and
isomorphisms. Each of the line bundles $\bar\L_\chi$ has a tautological section
$\sigma_\chi$.

For each pair $(\chi,\chi')\in G^*\times G^*$, the isomorphisms $\phi_j$
induce isomorphisms $$\phi_{\chi,\chi'}:\L_\chi\otimes\L_{\chi'}
\to \L_{\chi\chi'}\otimes\Bigotimes \M_i^{\epsi}.$$
Let $\tau_i\in H^0(V,\bar\M_i)$ be defined by $$
\tau_i=\sum_{\{\chi|(\ichi)\in S\}}\bar s_\ichi\sigma_\chi.$$
Define a section $\rho_{\chi,\chi'}$ of $\bar\L_\chi\otimes\bar\L_{\chi'}$ by
$$\rho_{\chi,\chi'}=\sigma_\chi\sigma_{\chi'}-
\bar\phi_{\chi,\chi'}^*(\sigma_{\chi\chi'}{\textstyle\prod}\tau_i^{\epsi}).$$
Then the zero locus of all the $\rho_{\chi,\chi'}$ is naturally a deformation
$\X\to T$ of
$X$ over $T$ (in particular $X$ can be naturally identified with the fibre of
$\X\to T$ over the closed point). This is proven in
\cite{Pa1} in the case where the deformation of
$Y$, $L_j$ and $M_i$ is the trivial one, but it is easy to see that the same
proof works in our generalized setting. The deformation $\X\to T$ so obtained
is
called the {\em natural deformation of $X$} associated to the given natural
deformation
of the reduced building data.}
\end{defn}
It is now clear why Galois deformations were called that way:
\begin{rem} Let $\X\to T$ be a deformation of $X$ induced by a
Galois deformation $(\Y,\ldots)$ of the reduced building data; $\X$ has a
canonical structure of $(G,I)$-cover of $\Y$, induced by the action of $G$ on
the total space of the line bundle $\L_\chi$ given by the character $\chi$.
\end{rem}
The restrictions to the category $\Art$ of the functors $\defnat $
and $\defgal$ satisfy Schlessinger's conditions for the existence of a
projective hull (see \cite{schl}); in fact, they can be described (as usual in
deformation
theory) in terms of tangent and obstruction spaces. If $F:\Art\to \Sets$ is a
contravariant functor, then we denote its tangent (resp.~obstruction) space by
$T^1(F)$ (resp.~$T^2(F)$), when this makes sense.
\setcounter{equation}{0}

\begin{lem}\label{tdefnat} There is a natural action of $G$ on $\defnat$,
whose invariant locus is $\defgal $; the decomposition of $T^l(\defnat )$
according to characters, for $l=1,2$, is the following:
\begin{eqnarray} &&T^l(\defgal )=T^l(\defnat )^{\rm
inv}=H^l(Y,T_Y(-\log{\textstyle\sum} D_i));\label{tgal}\\ &&T^l(\defnat
)^{\chi}=\Bigoplus_{i\in S_\chi}
H^{l-1}(Y,\O_Y(D_i)\otimes L_\chi^{-1})\qquad\hbox{for $\chi\ne1$;}
\end{eqnarray} where $S_\chi=\{i\in I|(\ichi)\in S\}$.
\end{lem}
\Pf An element $g\in G$ acts by $$(\Y,\M_i,\L_j,s_\ichi,\phi_j)
\mapsto
(\Y,\M_i,\L_j,\chi(g)s_\ichi,\phi_j).$$ It is clear that $\defgal$ is
contained in the invariant locus. It is not difficult to show the other
inclusion using the fact that the cover is totally ramified.

We now study separately tangent and obstructions spaces corresponding to
the different characters. For the trivial character, i.e. $\defgal$, the
functor is isomorphic to the deformation functor of the data $(Y,M_i,s_i)$;
(\ref{tgal}) is then well known (see \cite{We}).

Fix a nontrivial character $\chi$. Then the problem reduces to studying the
deformations of the zero section of a line bundle, given a
deformation of the base and of the bundle. The statement can then be
proven by applying the following lemma.
\qed
\begin{lem} Let $o\in B'\subset B\in\Art$ be schemes of length $1,n,n+1$
respectively for some $n$; for schemes, etc.~over $B$ denote the restriction to
$B'$ by a prime and the restriction to $o$ by ${}_o$. Let
$\Y\to B$ be a smooth projective morphism, $\L$ a line bundle on $\Y$; let
$s'$ be a section of $\L'$, such that $s'_o=0$. Then the obstruction to lifting
$s'$ to a section $s$ of $\L$ lies in $H^1(\Y_o,\L_o)$, and two liftings differ
by an element of $H^0(\Y_o,\L_o)$.
\end{lem}
\Pf Let $\{U_\alpha\}$ be an affine open cover of $Y=\Y_o$ such that $L$ is
trivial on each $U_\alpha$. Let $U_\albe$ be $U_\al\cap U_\be\subset U_\al$.
As $Y$ is smooth, we have that $\Y$ is covered by open subsets $V_\al$
isomorphic to
$U_\al\times B$, glued via $B$-isomorphisms $\phi_{\albe}:U_\albe\times B\to
U_\beal\times B$ satisfying the cocycle condition and restricting to the
identity over $o$. Let $g_{\albe}$ be transition functions for $\L$ with
respect to the open cover $V_\al$.

The section $s'$ can be described by functions $s'_\al$ on $U_\al\times B'$
such that, on $U_\albe\times B'$, $$
s'_\al=g'_\albe(s'_\be\circ\phi_\albe).$$
Extend $s'_\al$ arbitrarily to a function $s_\al$ on $U_\al\times B$;
any other extension is of the form $s_\al+\eps\si_\al$, where $\eps=0$ is an
equation of $B'$ in $B$ and $\si_\al$ is a function on $U_\al$ (as $\eps f=0$
for any function $f$ in the ideal of $o$ in $B$).
If an extension $s$ of $s'$ exists, then there must be functions $\si_\al$ on
$U_\al$ such that, on $U_\albe\times B$,
$$
s_\al+\eps\si_\al=g_\albe((s_\be+\eps\si_\be)\circ\phi_\albe).$$
Let $u_\albe=s_\al-g_\albe(s_\be\circ\phi_\albe)$. The restriction of
$u_\albe$ to $U_\albe\times B'$ is zero, hence $u_\albe$ is divisible by
$\eps$: let $u_\albe=\eps v_\albe$.
One can verify, using the fact that $s_o=0$, that $v_\albe$ is a cocycle
in $H^1(Y,\L_o)$: it is enough to check that $$
u_\albe+g_\albe(u_\bega\circ\phi_\albe)=u_\alga $$
on $U_{\albe\ga}$, for all triples $\al,\be,\ga$ of indices of the cover.
It is then immediate to verify that $v_\albe$ is the obstruction to lifting
$s'$ to $\Y$, and the statement about the difference of two liftings can be
proven in a similar way.
\qed

We now recall some properties of $\defor X$.  Let $\defg X:\Pointed\to \Sets$
be
the functor of deformations of $X$ together with the $G$ action.
\begin{lem} There is a natural action of $G$ on $\defor X$, whose
invariant locus is $\defg X$.
\end{lem}
\Pf Let $\X\to T$ be a deformation of $X$ over $(T,o)$; there is a given
isomorphism $i:X\to \X_o$. The action of an element $g\in G$ is given by
replacing $i$ with $i\circ \phi(g)$, where $\phi:G\to Aut(X)$ is the natural
action.

It is clear that if $G$ acts on a deformation $\X\to T$, then this belongs to
$\defg X$. The other implication follows from \cite{Ca2}, \S 7 or directly
from the fact that the automorphisms of $X$ and of its deformations are rigid.
\qed
Note that, as $X$ is of general type,
the $G$-action on $\defor X$ induces an action on the Kuranishi family
$\X\to B$ of $X$; the restriction of the Kuranishi family to the fixed
locus $B^G$ is universal for the functor $\defg X$ (compare (\cite{Pi},
(2.8) p.~19, \cite{Ca2}, \S 7).

Recall the following result from \cite{Pa1}.
\setcounter{equation}{0}

\begin{lem} Let $X$ be a smooth $(G,I)$-cover of $Y$ with building data
$(D_i,L_\chi)$. Then the decomposition according to characters of
$H^l(X,T_X)$ is as follows:
\begin{eqnarray} &&H^l(X,T_X)^{\rm inv}=H^l(T_Y(-\log \sum_{i\in I} D_i))\\
&&H^l(X,T_X)^\chi=H^l(T_Y(-\log \!\!\!\sum_{i\in S_\chi}\!\!\! D_i)\otimes
L_\chi^{-1})\qquad \hbox{\ if $\chi\ne1$}
\end{eqnarray}
where $S_\chi$ is the same as in lemma \rm{\ref{tdefnat}}.
\end{lem}
\Pf This follows immediately from proposition 4.1. in \cite{Pa1}.
\qed
\setcounter{equation}{0}
\begin{cor}\label{tdef} Assume that, for all $\chi\in G^*\setminus 1$, the
bundles $L_\chi$ and
$\Omega_Y^1\otimes L_\chi$ are ample. Then there are natural exact
sequences, for all $\chi\in G^*\setminus 1$:
\begin{eqnarray} &&0\to \textstyle\bigoplus\limits_{i\in S_\chi}
H^0(Y,\O(D_i)\otimes L_\chi^{-1})\to H^1(X,T_X)^\chi\to 0.\\ &&0\to
\textstyle\bigoplus\limits_{i\in S_\chi} H^1(Y,\O(D_i)\otimes L_\chi^{-1})\to
H^2(X,T_X)^\chi.
\end{eqnarray}
\end{cor}
\Pf Fix $\chi\ne 1$, let $D=\sum_{i\in S_\chi}D_i$, and consider the following
diagram of sheaves with exact rows and columns: $$
\begin{array}{ccccccccc}
&&&& 0 && 0 \\
&&&& \downarrow && \downarrow \\
&&&& \Bigoplus_{i\in S_\chi}^{\phantom{S_\chi}}\O_Y & = &  \Bigoplus_{i\in
S_\chi}\O_Y\\
&&&& \downarrow && \downarrow \\
0 & \longrightarrow & T_Y(-\log D) & \longrightarrow &\PP^*&\longrightarrow&
\Bigoplus_{i\in S_\chi}^{\phantom{S_\chi}}\O_Y(D_i)&\longrightarrow&0\\
&&\Vert&& \downarrow && \downarrow \\
0 & \longrightarrow & T_Y(-\log D) & \longrightarrow &T_Y&\longrightarrow&
\Bigoplus_{i\in S_\chi}^{\phantom{S_\chi}}\O_{D_i}(D_i)&\longrightarrow&0\\
&&&& \downarrow && \downarrow \\
&&&& 0 && 0
\end{array}$$
where $\PP$ is the prolongation bundle associated to the normal crossing
divisor $D$. By the previous lemma, it is enough to prove that the first two
cohomology groups of $\PP^*\otimes L_\chi^{-1}$ vanish; this follows from the
corresponding vanishing for $L_\chi^{-1}$ and $T_Y\otimes L_\chi^{-1}$, and
the latter is just Kodaira vanishing (it is here that one needs the assumption
$\dim Y\ge 2$).
\qed
The natural transformation of functors $\defnat\to \defor
X$ defined in \ref{trasnat} is equivariant with respect to the natural actions
of $G$ on these functors. Therefore, there is a commutative diagram
$$\begin{array}{ccc}
\defgal  & \longrightarrow & \defg X \\
\longdownarrow &    &\longdownarrow \\
\defnat  & \longrightarrow & \defor X

\end{array}
$$ where the vertical arrows are injections. The following theorem shows
that the horizontal arrows are smooth morphisms of functors when the branch
divisors are sufficiently ample.

This was proven in \cite{Pa1} under the hypothesis that $Y$ be rigid and
regular; in this case natural deformations are unobstructed, and it is
enough to check the surjectivity of the Kodaira-Spencer map. In the general
case one has to take into account the
obstructions  as well.
\begin{thm}\label{complete} Let $f:X\to Y$ be a totally ramified $(G,I)$-cover
with building
data $D_i$, $L_\chi$, such that $X$ and $Y$ are smooth of dimension $\ge 2$
and that $X$ is of general type. Assume that for all $\chi\in G^*\setminus 1$
the
bundles $L_\chi$ and $\Omega_Y^1\otimes L_\chi$ are ample. Then the natural
map of
functors (from $\Art$ to $\Sets$) $\defnat \to \defor X$ is smooth, and so is
the induced map
$\defgal $ and
$\defg X$.
\end{thm}
\Pf By a well-known criterion, smoothness of a natural transformation of
functors is implied by surjectivity of the induced map on tangent spaces,
and injectivity on obstruction spaces.

This is immediate by lemma \ref{tdefnat} and corollary \ref{tdef}, and by
the fact that the map between tangent (obstruction) spaces induced by the
map of functors is the natural one.
\qed

\section{Main theorem} In this section we will prove that the automorphism
group of an abelian cover is precisely the Galois group, provided that the
branch divisors are sufficiently ample and generic. The proof depends on the
construction of an explicit partial resolution of some singular covers, which
will be given in section 7.

Although the result is in some sense expected, the proof is rather involved
and the techniques applied are, we believe, of independent interest.

The following lemma is inspired by a similar result of McKernan (\cite{McK}).

\begin{lem} \label{mac} Let $\De$ be the unit disc in $\C$,
$\De^*=\De\setminus\{0\}$. Let $p:\X\to \De$ be a flat map, smooth over
$\De^*$, whose fibres are  integral
 projective varieties of non negative
Kodaira dimension. Assume we are given a section $\sigma$ of
$Aut_{\X/\De^*}$. If there exists a resolution of singularities $\eps:\tilde
\X\to\X$ such that each divisorial component of the exceptional locus has
Kodaira dimension $-\infty$, then $\sigma$ can be (uniquely) extended to a
section of $Bir_{\X/\De}$.
\end{lem}
\Pf The section $\sigma$ induces a birational map $\phi:\X\birto \X$ over
$\Delta$; the uniqueness of the extension follows from this. Let
$\tilde\phi:\tilde \X\birto\tilde \X$ be the induced birational map, and let
$\Gamma$ be a resolution of the closure of the graph of $\tilde\phi$; let
$p_1$, $p_2$ be the natural projections of $\Gamma$ on $\tilde\X$ (such that
$p_2=\tilde\phi\circ p_1$), and let $q_i=\eps\circ p_i$.

The strict transform $\X_0'$ of $\X_0$ in $\Gamma$ via $q_1^{-1}$ has
positive Kodaira dimension, hence it cannot be contracted by $p_2$, which is
a birational morphism with smooth image.  Therefore the restriction of $p_2$
to $\X_0'$ is birational (because $\X_0'$ is not contained in the exceptional
locus of $p_2$) onto some irreducible divisor $\X_0''$ in $\X$.

As $\X''_0$  is birational to $\X_0$ it cannot be of Kodaira dimension
$-\infty$; hence it is not contained in the exceptional locus of $\eps$.
Therefore $\eps(\X_0'')$ is a divisor contained in $\X_0'$, hence it is
$\X_0'$ by irreducibility, and the map $\eps:\X_0''\to\X_0$ is birational.

So the birational map $\phi$ can be extended to $\X_0$ by the birational map
$q_2\circ \left(q_{1|\X_0'}\right)^{-1}$.
\qed

\begin{lem} \label{trick}In the same hypotheses of lemma {\rm \ref{mac}},
assume moreover that there is a line bundle $L$ on $\X$, flat over $\De$,
whose restriction to $\X_t$ is very ample for all $t$, and such that
$h^0(\X_t,L_{|\X_t})$ is constant in $t$. If the action of $\sigma$ can be
lifted to an action on $L$, then $\sigma$ can be uniquely extended to a
section of $Aut_{\X/\De}$.
\end{lem}
\Pf Let $N$ be the rank of the vector bundle $p_*L$ on $\Delta$; choosing a
trivializing basis yields an embedding $\X\into\P^{N-1}\times \De$. The
automorphisms $\phi_t$ of $\X_t$ are restrictions to $\X_t$ of nondegenerate
projectivities of $\P^{N-1}$; their limit, as $t\to 0$, is a well-defined,
possibly degenerate projectivity $\phi_0$. This gives an extension of $\phi$
to an open set of $\X_0$; this must now be birational by the previous lemma,
which in turn implies that $\phi_0$ is nondegenerate (as $\X_0$ is not
contained in a hyperplane), and therefore that $\phi_0$ is a morphism.
Applying the same argument to $\phi^{-1}$ concludes the proof.
\qed
\begin{rem} The hypothesis that $\sigma$ acts on $L$ is obviously verified
if  $L_{|\X_t}$ is a pluricanonical bundle for all $t\ne 0$.
\end{rem}
\begin{prop}\label{prop.aut}
Let $p:\X\to \De$
be a flat family of integral projective varieties of general type, smooth
over $\De^*$. Assume that there is a line bundle $L$ on $\X$, flat over
$\De$, with $L_t:=L_{|\X_t}$ ample on $\X_t$, and $Aut(\X_t)$ acts on $L_t$
for $t\ne 0$. Assume
moreover that for any $m$-th root base change $\rho_m:\Delta\to \De$ the
pullback $\rho_m^*\X$ admits a resolution having only divisors of negative
Kodaira dimension in the exceptional locus. Then $Aut_{X/\De}$ is proper over
$\De$, and the cardinality of the fibre is an upper semi-continuous function.
\end{prop}
\Pf  After replacing $L$ with a suitable multiple and maybe shrinking $\Delta$,
we can assume that $L_t$ is very ample on $\X_t$, and that $h^0(\X_t,L_t)$
is constant in $t$. The map $Aut_{\X/\De}\to \De$ is obviously quasi-finite
(because the
fibres are of general type) and the fibres are reduced (because automorphism
groups are always reduced in char.~$0$). It is enough to prove that given a
map of a pointed curve $(C,P)$ to $\De$ and a lifting of the map to
$Aut_{\X/B}$ out of $P$, the lifting can be extended to $P$.

Via restriction to an open set we can assume that $C$ is the unit disc $\De$,
$P$ is the origin and $\De\to \De$ is the map $z\to z^m$; we can then apply
lemma \ref{trick} to conclude the proof.
\qed
\begin{cor}\label{gen.aut} Let $\X\to B$ be a smooth family of varieties
having ample canonical bundle. Then the scheme $Aut_{X/B}$ is proper over
$B$, and the cardinality of the fibre is an upper semi-continuous function.
\end{cor}
\Pf We can apply the previous proposition with $L=K_{\X/\De}$.
\qed
\begin{thm} \label{mainthm} Let $Y$ be a smooth projective variety, and $X$ a
smooth $(G,I)$-cover with ample canonical bundle, with covering data
$L_\chi$, $D_i$. Let $H=\O_Y(1)$ for some embedding of $Y$ in $\P^{N-2}$;
assume that the linear system $$|D_1-m_1NH|$$ is base-point-free. Assume also
that the $\Q$-divisor $$ M=K_Y-(m_1-1)NH+\sum_{i\in I}\frac{(m_i-1)
}{m_i}D_i$$ is ample on $Y$. Then, for a generic choice of $D_1$ in its linear
system, $X$
has automorphism group isomorphic to $G$.
\end{thm}
\newcommand{\iz}{1}
\Pf Let $d$ be the number of automorphisms of a generic cover with the given
covering data (cfr.~corollary \ref{gen.aut}). It is enough to show that $d\le
\#G$, the other inequality being obvious.

Let $H$ be as in the statement of the theorem, and let $\HH\subset |H|$ be
the (not necessarily complete) linear system giving the embedding; let
$H_1,\ldots,H_N$ be $N$ projectively independent divisors in $\HH$. Assume
that the $H_i$'s are generic, in particular that they are smooth and that
their union with all of the $D_i$'s has normal crossings. Let
$m=m_\iz$, $D=D_\iz$.

The strategy of the proof is the following: start from a generic cover $X$ of
$Y$, and construct a sequence of manifolds $X_1,\ldots,X_N$ and of subgroups
$G_k$ of $Aut(X_k)$ such that $$
\#Aut(X)\le \#G_1\le\ldots\le \#G_N\qquad{\rm and}\qquad G_N=G.$$ In fact,
$X_k$ will be a $(G,I)$-cover of $Y$ with covering data $D^{(k)},D_2,\ldots$,
$L_\chi^{(k)}$, where  $L_\chi^{(k)}=L_\chi-k\re^1_\chi H$ and $D^{(k)}$ is a
generic divisor in $|D-kmH|$ (recall that $\re^i_\chi$ was defined as the
unique integer $a$ satisfying $0\le a\le m_i-1$ and $\chi_{|H_i}=\psi_i^a$).
We let $G_k$ be the group of automorphisms of $X_k$ preserving the inverse
images of the curves $H_1,\ldots,H_k$ in $Y$.

We therefore want to prove the following:\begin{enumerate}
\item $\#Aut(X)\le \#G_1$;
\item $\#G_k\le \#G_{k+1}$;
\item $G_N=G$.
\end{enumerate}
\smallskip
\noindent {\sc First step:} $\#Aut(X)\le \#G_1$. Let $D^{(1)}$ be a generic
divisor in $|D-m H|$, and choose equations $f_1$, $g$ and $h_1$ for $H_1$,
$D$ and $D^{(1)}$ respectively. Define divisors $\D_i$ on $Y\times \C$ by
$\D_i=D_i\times \C$ for $i\ne \iz$, $\D_\iz=\{(1-t)f_1^mh_1+tg=0\}$; let
$\X^1$ be the corresponding abelian cover. $\X^1_0$ is a singular variety
(singular along the inverse image of the curve $H_1$ in $Y$), with
smooth normalization $X_1$ (see \cite{Pa1}, step 1 of normalization algorithm
of p.~203). Note that $X_1$ is of general type by the ampleness assumption on
$M$.

By proposition \ref{reslemma}, the family $\X^1$ and each $n$-th root base
change of $\X^1$ admit a resolution with only divisors of Kodaira dimension
$-\infty$ in the exceptional locus. Moreover, the pull-back of $(\#G)M$
restricts to the $\#G$-canonical bundle on the smooth fibres of $\X^1$ (cfr the
proof of prop.~4.2 in \cite{Pa1}, p.~208). Applying proposition \ref{prop.aut}
gives that $Aut_{\X^1/\C}$ is proper over $\C$, and hence that $\#Aut(X)\le
Aut(\X^1_0)$ (as we assumed $X$ to be generic). On the other hand it is clear
that each automorphism of $\X^1_0$ lifts to the normalization $X_1$, yielding
an automorphism which maps to itself the inverse image of the singular locus,
i.e., the inverse image of the curve $H_1$.

\smallskip
\noindent {\sc Second step:} $\#G_{k-1}\le \#G_{k}$. We use a similar
construction; let $X_{k-1}$ be as above, let $h_{k-1}$ be an equation of
$D^{(k-1)}$, $f_k$ an equation of $H_k$, and $h_k$ an equation of  $D^{(k)}$.
Define a $(G,I)$-cover $\X^k$ of $Y\times
\C$ branched over $D_i\times \C$ for $i\ne \iz$, and over
$\D_1^{(k)}=\{(1-t)f_k^mh_k+th_{k-1}=0\}$; $\X^k_0$ is singular along the
inverse image $C_k$ of $H_k$, and its normalization is $X_k$; again $X_k$ is
of general type.

Again by proposition \ref{reslemma} the family $\X^{(k)}$ and all its $n$-th
root base changes have a resolution with only uniruled components in the
exceptional locus; the same argument as before proves the result.

\smallskip\noindent {\sc Final step:} $G_N=G$. Let $\pi:X_N\to Y$ be the
covering map: $G_N$ is the group of automorphisms of $Y$ fixing the inverse
images of the curves $H_1,\ldots,H_N$. Every element of $G_N$ preserves
$\pi^*\left(\HH\right)$, hence induces an automorphism of $Y$; this
automorphism must be the identity as it induces the identity on $\HH$.
Therefore $G_N$ must coincide with $G$.
\qed
\begin{rem} In theorem {\rm \ref{mainthm}} we can replace the assumption that
the linear system $|D_1-m_1NH|$ be base point free by asking that for each
$i\in I$ $$ |D_i-m_iN_iH|$$ be base point free, with $N_i$ nonnegative
integers with sum $N${\rm ;} we then get that, for a generic choice of the
$D_i$'s such that $N_i\ne 0$, $Aut(X)=G$.
\end{rem}
\begin{ex}{\rm One might wonder whether it is always true that a generic
abe\-lian cover of general type has no ``extra automorphisms". Here is an easy
example where this is not the case. Consider a $\Z_3$-cover of
$\P^1$, branched over two pairs of distinct points, with opposite characters.
A generic such cover is a smooth genus $2$ curve, hence its automorphism
group cannot be $\Z_3$. }
\end{ex}
\begin{ex}{\rm Here is a slightly more complicated example of
extra  automorphisms, which works in any
dimension. Let $Y$ be a principally polarized abelian variety, and let $L$
be a principal polarization; assume that $L$ is symmetric, i.e. invariant
under the natural involution $\sigma(y)=-y$ on $Y$. The sections of
$L^{\otimes 2}$ are all symmetric, and the associated linear
system  has no base points. Let $G=\Z_2^s$, with the canonical basis
$e_1,\ldots,e_s$. Choose $I=\{1,\ldots,s\}$, and let $H_i$ be the subgroup
generated by $e_i$, for $i=1,\ldots,s$.
\par
The equations for the reduced building data become $L_j^{\otimes
2}=\O_Y(D_j)$; we choose the solution $L_j=L$, $M_i=L^{\otimes 2}$ for all
$i,j$. We are in fact constructing a fibred product of double covers. Choose
the $D_i$'s to be generic divisors in the linear system $|L^{\otimes 2}|$.
Each of them must be symmetric; this implies that the involution $\sigma$ can
be lifted to an involution of $X$, which is an automorphism not contained in
the Galois group of the cover.
\par
Note
that in this case the total branch divisor can become arbitrarily large, still
all
$(G,I)$-covers have an automorphism group bigger than $G$. }
\end{ex}

\section{Moduli spaces of abelian covers and global constructions}

In this section we will explicitly
construct a coarse moduli space for abelian covers of a smooth variety $Y$
and a complete space of natural deformations. Although some of the material in
this section is implicit in
\cite{Pa1}, we find it important to state it in a precise and explicit way.
In particular we will apply theorem \ref{complete} to construct (under
suitable ampleness assumptions) a family of natural deformations which maps
dominantly to the moduli (theorem 5.12).

Let $Y$ be a smooth, projective variety, $G$ an abelian group, $I$ a subset
of $I_G$. A {\em family
of smooth $(G,I)$-covers} of $Y$ over a base scheme $T$ is a smooth, proper
map
$\X\to T$ and an action of $G$ on $\X$ compatible with the projection on $T$,
together with a $T$-isomorphism of the quotient $\X/G$ with $Y\times T$,
such that for each $t\in T$ the induced cover $\X_t\to Y$ is a
$(G,I)$-cover. Two families over
$T$ are {\em {\rm (}strictly\/{\rm )} isomorphic} if there is a
$G$-equivariant isomorphism inducing on the quotient $Y\times T$ the identity
map.

A (coarse) moduli space $\ZZ$ for smooth $(G,I)$-covers of $Y$ is a
scheme structure on the set of smooth $(G,I)$-covers modulo isomorphisms,
such that for any family of $(G,I)$-covers of $Y$ with base $T$ the induced
map $T\to
\ZZ$ is a morphism.

\begin{thm} There is a coarse moduli space of $(G,I)$-covers of $Y$, which is
a Zariski open set $\ZZ=\ZZ(Y,G,I)$ in the closed subvariety of
$$\prod_{\chi\in G^*\setminus 1}Pic(Y)\times \prod_{i\in I} \hilbd(Y)$$ of
all the $(L_\chi,D_i)$ satifying the relations {\rm (\ref{bdata})}. The open
set $\ZZ$ is the set of $(L_\chi,D_i)$'s which satisfy the additional
conditions:
\begin{enumerate}
\item each $D_i$ is smooth and the union of the $D_i$'s is a divisor with
normal crossings;
\item whenever $D_{i_1}, \ldots, D_{i_k}$ meet, the natural map
$H_{i_1}\oplus\cdots\oplus H_{i_k}\to G$ is injective.
\end{enumerate}
\end{thm}
\Pf The set $\ZZ$ parametrizes the smooth abelian covers of $Y$ by
\cite{Pa1}, theorem 2.1. The fact that the induced maps from a family of
abelian covers to $\ZZ$ are morphisms follows from the corresponding property
of the Hilbert schemes and Picard groups.
\qed

Proposition 2.1 of \cite{Pa1} implies:
\begin{rem} For any basis $\chi_1,\ldots,\chi_s$ of $G^*$, the natural map
$$\ZZ\to\prod_{j=1}^s Pic(Y)\times\prod_{i\in I} \hilbd(Y)$$ induced by
projection is an isomorphism with its image.
\end{rem}

$\ZZ$ decomposes as the disjoint union of infinitely many quasiprojective
varieties $Z(\xi_i,\eta_\chi)=Z(\xi_i,\eta_\chi)(Y,G,I)$, where $\eta_\chi$,
$\xi_i$ are the Chern classes of $L_\chi$ and $\O(D_i)$, respectively. We now
give an explicit description of $Z(\xi_i,\eta_\chi)$ under the assumption
that the
$\xi_i$'s are sufficiently ample.
\setcounter{equation}{0}
\begin{prop}\label{moduli} Let $\xi_i$, $\eta_\chi$ be cohomology classes
satisfying the following relations {\rm(}compare {\rm(\ref{bdata})):}
\begin{equation}
\label{Chbdata}\eta_\chi+\eta_{\chi'}=\eta_{\chi\chi'}+\sum_{i\in I}
\epsi\xi_i\qquad\quad \hbox{for
all $\chi,\chi'\in G^*\setminus 1$}.
\end{equation} Assume moreover that $\xi_i-c_1(K_Y)$ is the class of an ample
line bundle for all $i\in I$. Then $Z(\xi_i,\eta_\chi)$ is an open set in a
smooth fibration (with fibre a product of projective spaces) over an abelian
variety $A(\xi_i,\eta_\chi)$ isogenous to $Pic^0(Y)^{\#I}$. $Z(\xi_i,
\eta_\chi)$ is nonempty iff
there are smooth effective divisors $D_i$, with  $c_1(D_i)=\xi_i$, such that
their union has normal crossings.
\end{prop}
\Pf Let $A=A(\xi_i,\eta_\chi)\subset\prod_{i\in
I}Pic^{\xi_i}(Y)\times \prod_{\chi\in G^*\setminus 1}Pic^\chi(Y)$ be the image
of $Z(\xi_i,\eta_\chi)$; by
equations (\ref{rbdata}) the natural map $A\to \prod_{i\in I}Pic^{\xi_i}(Y)$
is a finite \'etale cover of degree $(2q)^{\#G}$, where $q$ is the
irregularity of $Y$. So each connected component of $A$ is an abelian
variety, isogenous to $Pic^0(Y)^{\#H}$. The fact that $A$ is connected is a
consequence of the covering being totally ramified. In fact, choose a basis
$\chi_1,\ldots,\chi_s$ of $G^*$, and consider the diagram
$$
\begin{array}{ccc} A&\longrightarrow&\prod_{i\in I}Pic^{\xi_i}(Y)\\
\longdownarrow& &\longdownarrow\\
\prod_{j=1}^sPic^{\eta_j}(Y)&\longrightarrow&\prod_{j=1}^sPic^{\ord{\chi_j}
\eta_j}\\
\end{array}$$ with maps given by $$
\begin{array}{ccc} (M_i,L_j)&\mapsto&(M_i)\\
\downarrow& & \downarrow\\ (L_j)&\mapsto&(L_j^{\otimes n_j}=\otimes M_i^
{\reb^i_{j}}).\\
\end{array}$$ The diagram is a fibre product of (connected) abelian varieties;
to prove that $A$ is connected is equivalent to proving that
$\pi_1(\prod_{i\in I}Pic^{\xi_i}(Y))$ surjects on
$$\pi_1(\prod_{j=1}^sPic^{\ord{\chi_j}\eta_j})/
\pi_1(\prod_{j=1}^sPic^{\eta_j}(Y));$$ this is in turn equivalent to proving
that $G^*$ injects in $\oplus_{i\in I}H_i^*$, which follows by dualizing
from assumption \ref{totram}.

Let $\PP_i$ on $A\times Y$ be the pullback of the Poincar\'e line bundles
from $Pic^{\xi_i}(Y)\times Y$; the pushforward of $\PP_i$ to $A$ is a vector
bundle $E_i$ because of the ampleness condition (the rank of $E_i$ can be
computed by Riemann-Roch).  The moduli space $Z(\xi_i,\eta_\chi)$ is an open
set of the fibred product of the $\P(E_i)$.
\qed

\begin{rem} {\rm If $q(Y)$ is not zero, then the components
$Z(\xi_i,\eta_\chi)$ are uniruled, but not unirational.}
\end{rem}

\begin{rem}{\rm In general $Z(\xi_i,\eta_\chi)$ is a coarse
but not a fine moduli space, i.e., it does not carry a universal family.
Keeping the notation of proposition \ref{moduli}, let $\VV$ be the total
space of the fibred product of the $E_i$'s, and let $\VV^o$ the inverse image
of
$Z(\xi_i,\eta_\chi)$; we have a natural abelian cover of $Y\times \VV^o$,
which is a complete family of smooth covers of $Y$ with the given data.}
\end{rem}

There is a natural action of $Aut(Y)$ on the moduli space of $(G,I)$-covers
$\ZZ$, given by
$$\phi(D_i,L_\chi)=(\phi(D_i),(\phi^{-1})^*L_\chi)\qquad\qquad\hbox{for
$\phi\in Aut(Y)$.}$$ The automorphism group of $G$ acts naturally on $G^*$
(by $\autg(\chi)=\chi\circ \autg^{-1}$) and on $I_G$ (by $\autg
(H,\psi)=(\autg(H),\psi\circ \autg^{-1})$); given a subset $I$ of $I_G$, let
$Aut_I(G)$ be the set of automorphisms of $G$ preserving $I$. There is a
natural action of $Aut_I(G)$ on $\ZZ$, induced by the natural action of this
group on the indexing sets $G^*\setminus 1$ and $I$.

\begin{prop} If the classes $\xi_i$'s are ample enough {\rm (}so that theorem
{\rm \ref{mainthm}} applies to some cover in $Z(\xi_i,\eta_j)${\rm),} then
the quotient
of
$Z(\xi_i,\eta_j)$ by the natural action of
$Aut(Y)\times Aut_I(G)$ maps birationally to its image in the moduli of
manifolds with ample canonical class.
\end{prop}
\Pf That the natural map to the moduli factors via this action is clear.
Viceversa, given a generic cover $X$ in $Z(\xi_i,\eta_j)$, by theorem
\ref{mainthm} its automorphism group is isomorphic to $G$; so it can be
identified uniquely as a $(G,I)$-cover up to isomorphisms of $G$ and of $Y$.
\qed
\begin{defn}{\rm Let
$\Y\to T$ be a deformation of $Y$ over a simply connected pointed analytic
space $(T,o)$. As $T$ is simply connected, the cohomology of every fibre
$\Y_t$ is canonically isomorphic with that of $Y$. Then the varieties
$Z(\xi_i,\eta_\chi)(\Y_t,G,I)$ (resp.\ $A(\xi_i,\eta_\chi)(\Y_t,G,I)$) for
$t\in T$ glue to a global variety
$\ZZ_T(\xi_i,\eta_\chi)=\ZZ_T(\xi_i,\eta_\chi)(\Y,G,I)$ (resp.\
$\A_T(\xi_i,\eta_\chi)$), surjecting on the
locus on $T$ where the classes
$\xi_i$ {\rm(}and hence also the $\eta_\chi${\rm)} stay of type $(1,1)$. The
global varieties are constructed by
replacing the Hilbert and Picard schemes in the construction of
$Z(\xi_i,\eta_\chi)$ and $A(\xi_i,\eta_\chi)$ with their relative versions.
The previous results can all be extended to this relative setting.
}
\end{defn}

For each smooth $(G,I)$-cover $f:X\to Y$, the natural deformations of
the reduced building data such that the induced deformations of
$(Y,L_j,M_i)$ is trivial are parametrized naturally by $\prod_{(\ichi)\in
S}H^0(Y,M_i\otimes L_\chi^{-1})$, as in \S 5 of \cite{Pa1}.
\begin{thm} Let $\Y\to T$ be a deformation of $Y$ over a germ $(T,o)$, and
assume that the $\xi_i$'s stay of type $(1,1)$ on $T$. Then there is a
quasiprojective morphism  $\W_T(\xi_i,\eta_\chi)\to \A_T(\xi_i,\eta_\chi)$
whose fibre over a point parametrizing line bundles $(L_j,M_i)$ on $\Y_t$
is canonically isomorphic to $\prod_{(\ichi)\in
S}H^0(Y,M_i\otimes L_\chi^{-1})$.
\end{thm}
\Pf The theorem follows, by taking suitable fibre products, from the
following two lemmas.
\qed
\begin{lem} Let $Y$ be a smooth projective variety, and $\xi\in NS(Y)$. Then
there exists a morphism of schemes $\pi:W^\xi(Y)\to Pic^\xi(Y)$ such that the
fibre over a point $[L]$ is naturally isomorphic to the vector space
$H^0(Y,L)$. For any choice of the Poincar\'e line bundle $\PP$ on $Y\times
Pic^\xi(Y)$, there exists such a $W^\xi(Y)$ with the property that the line
bundle $\pi^*\PP$ on
$Y\times W^\xi(Y)$ has a tautological section.
\end{lem}
Let $\PP$ be the Poincar\'e line bundle on $Y\times Pic^\xi(Y)$, and let
$p:Y\times Pic^\xi(Y)\to Pic^\xi(Y)$ and $q:Y\times Pic^\xi(Y)\to Y$ be the
projections; if
$p_*(\PP)$ is a vector bundle, it is enough to take $W$ to be the
total space of this vector bundle.

It is also clear that if $\xi-c_1(K_Y)$ is an ample class, then $p_*(\PP)$
is indeed a vector bundle. For the general case, let $A$ be a line bundle on
$Y$ such that $c_1(A)+\xi-c_1(K_Y)$ is ample, and such that there exists an
$s\in H^0(Y,A)$ defining an effective, smooth divisor $D$. Let $\pi:V\to
Pic^\xi(Y)$ be
the total space of the vector bundle $p_*(\PP\otimes q^*A)$, and let
$\sigma:\O_{Y\times V}\to
\pi^*(\PP\otimes q^*A)$ be the tautological section.
For every $y\in D$, let $\sigma_y$ be the induced section of $\pi^*(\PP\otimes
q^*A)|_{\{y\}\times V}$; let $W_y\subset V$ be the divisor defined by
$\sigma_y$. Let $W=W^\xi(Y)$ be the intersection of all $W_y$'s for $y\in D$:
then
$\sigma/s$ is regular on $W$, and defines the required tautological
section.
\qed
\begin{lem}{Let $\Y\to T$ be a deformation of $Y$ over
a germ of analytic space $T$, and assume that $\xi$ stays of type $(1,1)$
over $T$. Then, after maybe replacing $T$ with a Zariski-open subset, the
spaces $W^\xi(\Y_t)$ glue together to a quasiprojective morphism
$W^\xi_T(\Y)\to Pic^\xi_T(\Y)$.}
\end{lem}
\Pf After possibly restricting $T$, we can extend $A$ to a line bundle $\A$
over
$\Y$, and $s$ to a section of
$\A$. The rest of the proof remains the same, using the fact that the
relative Picard scheme exists and carries a Poincar\'e line bundle.
\qed
We now want to describe explicitly
$W^\xi(Y)$ in the case $\xi=0$, which we will use repeatedly later.
\begin{rem}{\rm For any deformation $\Y\to T$ over a germ of analytic space,
$W^0_T(\Y)$ is naturally isomorphic to the union in $Pic^0_T(\Y)\times \C$
of
$j(T)\times\C$ and
$Pic^0_T(\Y)\times
\{0\}$, where $j:T\to Pic^0_T(\Y)$ is the zero section.}
\end{rem}
In particular $W^0(Y)$ is reducible when $q(Y)\ne 0$; this reflects the fact
that the deformations, as pair (line bundle, section), of
$(\O_Y,0)$ are obstructed; one can either deform the line bundle or the
section, but not both at the same time. This remark will be used to construct
examples of manifolds lying in several components of the moduli in section 6.

\begin{thm} {\rm (i)} Let $Y$ be a smooth projective variety, and let $X\to Y$
be a smooth $(G,I)$-cover such that theorem {\rm \ref{complete}} holds. Then
there exists a pointed analytic space $(\W,w)$ and a natural deformation of
the reduced building data of $X$ over $\W$ such that the induced map of germs
from $(\W,w)$ to the Kuranishi family of $X$ {\rm(}defined as in {\rm
3.3}{\rm)} is surjective.
\par\noindent
{\rm (ii)} One can choose $\W$ to be a quasi-projective scheme, and then the
induced rational map from $\W$ to the moduli of manifolds with ample canonical
class is dominant onto each component of the moduli containing $[X]$.
\end{thm}
\Pf
(i) Let $\Y\to T$ be the
restriction of the Kuranishi family of $Y$ to the locus where all the
$\xi_i$'s stay of type $(1,1)$. Let $\W=\W_T(\xi_i,\eta_\chi)$,
$\Y_\W=\Y\times_T\W$. Over $\Y_\W$ there are tautological line bundles
$\L_j$, $\M_i$ and tautological sections $s_{\ichi}$ of
$\M_i\otimes\L_\chi^{-1}$ (where $\L_\chi$ is defined as in \ref{natdef});
moreover
$\L_j^{\otimes n_j}$ is isomorphic to $\Bigotimes \M_i^{\reb^i_j}$. $\W$
parametrizes data $(\Y_t,L_j,M_i,s_{\ichi})$ such that $t\in
T$, $L_j$ and $M_i$ are line bundles on $\Y_t$ satisfying (\ref{rbdata}) and
having Chern classes $\eta_j,\xi_i$, and $s_\ichi$ are sections of
$L_\chi\otimes M_i^{-1}$.
Let $w\in \W$ be a point corresponding to the reduced building data of $X$:
that is, assume that $w$ corresponds to the data $(\Y_o,L_j,M_i,s_{\ichi})$,
where
$s_\ichi=0$ for all
$\chi\ne 1$, $o$ is the chosen point in $T$, and the sections $s_{i,0}$
define divisors $D_i$ such that $(L_j,D_i)$ are the reduced building data of
$X$.

Choose arbitrarily isomorphisms $\Phi_j:\L_j^{\otimes n_j}\to \Bigotimes
\M_i^{\otimes\reb^i_j}$, extending the isomorphism over $w$ induced by
multiplication in $\O_X$.

By theorem \ref{complete}, together with Artin's results on approximation of
analytic mappings (see \cite{Ar}), it is enough to show that every natural
deformation of the reduced building data of $X$ over a germ of analytic
space can be obtained as pullback from $(\W,w)$.

It is clear that all small deformations of the data
$(Y,L_j,M_i,s_\ichi)$ can be obtained as pullback from $W$. So it is enough
to prove that, up to isomorphism of natural deformations, we can choose the
$\phi_j$'s arbitrarily. This is proven in lemma \ref{lautnonconta}.

\noindent (ii) Start by noting that one can construct a deformation $\Y\to B$
of $Y$ over a pointed quasi-projective variety $(B,o)$, such that the germ of
$B$ at $o$ maps surjectively to the locus in the Kuranishi family of $Y$ where
the classes $\xi_i$'s stay of type $(1,1)$. In fact, choose any $\chi\in
G^*\setminus 1$, and let $L$ be a sufficiently big multiple of $L_\chi$;
assume in particular that $L$ is very ample and that all its higher cohomology
groups vanish. Let $N=\dim H^0(Y,L)-1$; choosing a basis of $  H^0(Y,L)$ gives
an embedding of $Y$ in $\P^N$. Take the union of the irreducible components of
the Hilbert scheme of $\P^N$ containing $b=[Y]$, and consider inside it the
open locus $B'$ of points corresponding to smooth subvarieties. Then the
natural map from the germ of $B'$ at $b$ to the Kuranishi family of $X$
surjects on the locus where $\eta_\chi$ stays of type $(1,1)$. Let $B$ be the
closed subscheme of $B'$ where also the classes $\xi_i$ stay of type $(1,1)$.

Let $\Y\to B$ be the universal family; by replacing $B$ with an \'etale open
subset we can assume that $\Y\to B$ has a section. Then (compare for instance
\cite{Mu}, p.~20) there exists a global projective morphism $\A\to B$ and line
bundles $\M_i$, $\L_j$ on $\Y\times_B\A$, such that $\A_b$ parametrizes line
bundles $(M_i,L_j)$ on $\Y_b$ such that firstly, they satisfy the usual
compatibility conditions, and secondly, the Chern classes of $(M_i,L_j)$ lie
in the orbit of $(\xi_i,\eta_j)$ via the monodromy action of $\pi_1(B,b)$.

Mimicking the proof in the germ case, and replacing $B$ by an \'etale open
subset if necessary, one can find a quasi-projective morphism $\W\to \A$ whose
fibre over a point corresponding to line bundles $(M_i,L_j)$ on $\Y_b$ is
isomorphic to $\prod H^0(\Y_b,M_i\otimes L_\chi^{-1})$ for $(\ichi)\in S$,
together with tautological sections $\sigma_{\ichi}$ of the pullbacks to
$\Y\times_B\W$ of
$\M_i\otimes\L_\chi^{-1}$.

Let $w\in \W$ be a point corresponding to the building data of $X$ as before.
Again (possibly passing to an \'etale open subset) one can extend the
multiplications isomorphisms
$\phi_j$ to isomorphisms
$\Phi_j:\L_j^{\otimes n_j}\to \Bigotimes \M_i^{\otimes\reb^i_j}$.

Putting everything together, we have a natural deformation of the building
data of $X$ over $(\W,w)$; this induces by (3.3) a rational map to the moduli
of manifolds with ample canonical class, which is a morphism on the open
subset of $\W$ where the natural deformation of $X$ is smooth. Applying the
same methods as in (i) implies that the map from $\W$ to the moduli is
dominant on each irreducible component containing $[X]$.
\qed
\begin{lem}\label{lautnonconta}
Let $T$ be a germ of analytic space. For any
$(\Y,\M_i,\L_j,s_{\ichi},\phi_j)\in\defnat(T)$, and for any other admissible
choice of isomorphisms $\phi_j':\L_j^{\otimes n_j}\to\bigotimes M_i^{\otimes
\re^i_j}$, there exist sections $s'_\ichi$ such that $
(\Y,\M_i,\L_j,s_{\ichi},\phi_j)$ is isomorphic to
$(\Y,\M_i,\L_j,s'_{\ichi},\phi'_j)$.
\end{lem}
\Pf
It is enough to show that there are automorphisms $\psi_i$ of $\M_i$ such
that the composition $(\bigotimes \psi_i^{\otimes \reb^i_j})\circ \phi_j$
equals $\phi_j'$; in fact in this case one can choose
$s'_\ichi=\psi_i^*(s_\ichi)$, for all $(\ichi)\in S$. 

As both $\phi_j$ and $\phi_j'$ are isomorphisms, $\phi_j=f_j\phi_j'$, where
$f_j$ is an invertible function on $\Y$ restricting to $1$ on the central
fibre.
Finding the $\psi_i$'s is equivalent to finding functions $g_i$'s on $\Y$ such
that $g_i$ restricts to $1$ on the central fibre and $f_j=\prod
g_i^{\reb^i_j}$,
for all $j=1,\ldots,s$. The existence of such $g_i$'s follows from the fact
that the matrix
$\re^i_j$ has rank equal to $s$, which in turn is implied by the cover
being totally ramified (see lemma 2.1).
\qed

\setcounter{equation}{0}
\begin{rem}{\rm There is a natural action of $(C^*)^{\#I}$ on the functor of
natural deformations, which is the identity on $(\Y,\L_j,\M_i,\phi_j)$ and
acts on
$\sigma_{i,\chi}$
by
\begin{equation}\label{C*action} (\lambda_i)_{i\in
I}(\sigma_{j,\chi})=\prod_{i\in I}\lambda_i^{\delta_{ij}m_i-\re^i_\chi}\cdot
\sigma_{j,\chi};\end{equation} This action has the property that the induced
flat maps $\X\to T$ are invariant under it; in particular the
natural map from $\W(\xi_i,\eta_j)$ to the moduli factors
through the corresponding action.
}
\end{rem}

\section{Applications to moduli}

In this chapter we want to apply the results on deformation theory together
with theorem \ref{mainthm} to study the generic automorphism group of some
components of the moduli spaces of manifolds with ample canonical class,
components containing suitable abelian covers with sufficiently ample branch
divisors.

To begin with, we study the case of simple cyclic covers (i.e., those for
which the Galois group
$G$ is cyclic and there is only one irreducible branch divisor).
\begin{prop}\label{cyclic} Let $f:X\to Y$ be a smooth simple cyclic cover,
with Galois group $\Z_m$, and reduced building data $D$ and $L$ (where $D$ is
a smooth divisor and $L$ is a line bundle satisfying $mL\equiv D$). Assume
that $D$ is sufficiently ample. Let $M$ be an irreducible component of the
moduli space of surfaces of general type containing $X$. Then $G_M$ is
trivial if $m\ge 3$, and $G_M=G$ if $m=2$.
\end{prop}
\Pf  In case $m=2$, it is easy to check that $H^i(X,T_X)$ is $G$-invariant
for $i=1,2$; hence the natural map $\defgal
\to \defnat$ is surjective,
and all deformations are Galois. By theorem \ref{mainthm}, $Aut(X)=G$ for a
generic choice of $D$ in its linear system.

If $m\ge 3$, assume without loss of generality that $D$ is generic in its
linear system. Let $(G,\chi)$ be the element of $I_G$ corresponding to the
only nonempty branch divisor. Then the natural deformations of $X$ such that
$Y$ and $\O(D)$ are fixed are parametrized by $$
\Bigoplus_{i=0}^{m-2} H^0(Y,L^{-i}(D))=\Bigoplus_{i=0}^{m-2} H^0(Y,L^{m-i});$$
in particular they are unobstructed. Moreover, given any nontrivial element
$g$ of the Galois group $G$, it acts on the (necessarily nonzero) summand
$H^0(Y,L^{m-1})$ as multiplication by $\chi(g)$, hence nontrivially;
therefore  $g$ does not extend to the generic deformation. By genericity
however $Aut(X)=G$, hence by semicontinuity of the automorphism group the
proof is complete.
\qed
Hence, to get nontrivial examples, and to prove the results on the
moduli claimed in the introduction, it is necessary to study more general
abelian covers.

\begin{constr}\label{construction} {\rm Let $s$ be an integer $\ge2$. Let
$d_1,\ldots,d_s$ be integers $\ge 2$, such that $d_i|d_{i+1}$ for $i\le s-1$;
let
$d_0=d_s$, and define integers $b_i$ by requiring that $b_id_i=d_0$, for all
$i=1,\ldots,s$. Let $G=\Z_{d_1}\times\cdots\times \Z_{d_s}$, and let
$e_1,\ldots,e_s$ be the canonical basis of $G$; let $\chi_1,\ldots,\chi_s$ be
the dual basis of $G^*$.
\par Let $e_0:=-(e_1+\ldots+e_s)$, and let $H_i$ be the subgroup generated by
$e_i$; for $i=0,\ldots,s$, let $\psi_i\in H_i^*$ be the unique character such
that $\psi_i(e_i)=\zeta_{d_i}$; note that, for each
$j=1,\ldots,s$ and $i\ne 0$ we have
$\re^i_j=\delta_{ij}$, while $\re^0_j=b_j(d_j-1)$. Moreover $\ord{e_i}=d_i$,
for
$i=0,\ldots,s$.  Let $I=\{0,\ldots,s\}$; identify $I$ with a subset of $I_G$
via
$i\mapsto (H_i,\psi_i)$.
\par Fix a smooth projective variety $Y$ of dimension $d$, and assume that
$s\ge d\ge 2$. Let $f:X\to Y$ be a $(G,I)$-cover of $Y$, with branch divisors
$D_i$. Equations (\ref{rbdata}) become $$
L_j^{\otimes d_j}=M_j\otimes M_0^{\otimes (d_j-1)}$$
for all $j=1,\ldots,s$, hence they can be solved by letting
$L_j=M_0\otimes F_j$, $M_j=M_0\otimes F_j^{\otimes d_j}$, for all $j=1,\ldots
,s$.
\par We compute explicitly $L_\chi$ for $\chi\in G^*$, using equation
(\ref{chidarbd}). Let $\chi\in G^*$, and write
$\chi=\chi_1^{\alpha_1}\cdots\chi_s^{\alpha_s}$, with $0\le
\alpha_j< d_j$. One gets $$
L_{\chi}=\Bigotimes_{j=1}^s F_j^{\otimes \alpha_j}\otimes M_0^{\otimes
N_\chi},$$ where
$N_\chi=-[(-\al_1b_1-\ldots-\al_sb_s)/d_0]$. In particular $N_\chi$ is an
integer $\ge 0$; $N_\chi=0$ if and only if $\chi=0$, $N_\chi=1$ if and only
if $\sum(\al_ib_i)\le d_0$.

In the following we will always assume that $c_1(F_j)=0$, for $j=1,\ldots,s$;
let $\xi=c_1(M_0)$.
Assume also that $X$ is a smooth $(G,I)$-cover, that is that the divisors
$D_i$ are smooth and their union has normal crossings.\par
\setcounter{equation}{0} In the surface case, one can compute the Chern
invariants of the cover $X$:
\begin{eqnarray*}
K_X^2/\#G&=&\left(K_Y+(s-(d_0^{-1}+\ldots+d_s^{-1}))\xi\right)^2\\
c_2(X)/\#G&=&c_2(Y)-((s+1)-(d_0^{-1}+\ldots+d_s^{-1}))\xi K_Y+\\
&&\left({s+2\choose 2}+
\sum_{i=0}^sd_i^{-1}+\sum_{0\le i<j\le s}d_i^{-1}d_j^{-1}\right)\xi^2.
\end{eqnarray*} The first equality follows from \cite{Pa1}, proposition 4.2;
the second from the additivity of the Euler characteristic, by decomposing
$Y$ in locally closed subsets according to whether a point lies in $2$, $1$
or no branch divisor. Note that no other possibilities can occur, as we
assume that the union of the branch divisors has normal crossings. The second
equality could also be derived by Noether's formula and proposition 4.2 in
\cite{Pa1}. }
\end{constr}
\begin{lem}\label{keylemma} Let $f:X\to Y$ be a $(G,I)$-cover as in
construction \ref{construction}.  Assume that $q(Y)$ is
nonzero, that
$\xi\in NS(Y)$ is sufficiently ample, that $F_j=\O_Y$ (for $j=1,\ldots,s$),
and that $D_i\in |M_i|=|M_0|$ is generic {\rm(}for
$i=0,\ldots,s${\rm)}.   Then, for each
$k=0,\ldots, s$, there exists a component $M_k$ of the moduli of
manifolds with ample canonical class, containing $X$, such that the generic
automorphism group $G_{M_k}\subset G$ is equal to
$G_k=\Z_{d_{1}}\times\ldots\times\Z_{d_k}$.
\end{lem}
\Pf By assumption $X$ has ample canonical class,
$Aut(X)=G$ and the natural deformations of $X$ are complete. Assume first that
$Y$ is rigid. Let $\chi\in G^*$ be such that
$N_\chi=1$, and let $(\ichi)\in S$; these are the only values of
$\ichi$ (with $\chi$ nontrivial) for which $M_i\otimes L_\chi^{-1}$ can have
sections, i.e.~can contribute to non-Galois deformations. In fact
$c_1(M_i\otimes L_\chi^{-1})=0$, hence it has sections if and only if it is
trivial (compare remark 5.11). The condition that the line bundle $M_i\otimes
L_\chi^{-1}$ be trivial can be expressed, in terms of the $F_j$'s, as
\begin{equation}\label{trivial}
\sum_j \al_j F_j= d_i F_i.\end{equation}
Let $T_k\subset Pic^0(Y)^s$ be the locus where $F_i=0$ for all $i>k$. Note
that $F_i=0$ for all $i>k$ implies that $M_i=M_0$ for all $i>k$, and that
$M_i\otimes L_\chi^{-1}$ is trivial for any $(\ichi)$ such that $N_\chi=1$,
$i>k$ and $\chi$ restricted to $G_k$ is trivial.

For a
generic choice of
$(F_j)\in T_k$, the line bundles $M_i\otimes L_\chi^{-1}$ are nontrivial for
each $\chi$ such that $\chi_{|G_k}\ne 1$; in fact, for any such $\chi$ there
exists $j_0\le k$ such that $\al_{j_0}> 0$, hence the coefficient of $F_{j_0}$
in (\ref{trivial}) is nonzero (being either $\al_{j_0}>0$ or
$\al_{j_0}-d_{j_0}<0$).

On the other hand, for each $j>k$, one has
$(0,\chi_j)\in S$ and $M_0\otimes L_j^{-1}$ is trivial (in fact one has to
exclude here the case where $d_0$ is equal to $2$, and hence all $d_i$'s are;
this case needs a slightly different analysis, see below). Hence for every
$g\in G\setminus G_k$, and for any $(G,I)$-cover with building data in $T_k$,
there are natural deformations of the cover to which the action of $g$ does
not extend.

Therefore the $(G,I)$-covers whose building data are in $T_k$, together with
their
natural deformations such that $s_{\ichi}=0$ for all $\chi$ acting
nontrivially on $G_k$, form an irreducible component of the Kuranishi family
of $X$; in fact they are parametrized by an irreducible variety, and at some
point they are complete (at least at all points corresponding to
$(G,I)$-covers with a generic choice of the $F_j$'s for $j\le k$). The generic
element of this component has therefore automorphism group $G_k$.

In the case where $d_0=2$, $(0,\chi_j)\notin S$; however, if $k\ne s-1$, we
can consider $M_{j'}\otimes L_j^{-1}$ instead of $M_0\otimes L_j^{-1}$, where
$j'$ is any index $>k$ and different from $j$. If $k=s-1$, let
$\chi=\chi_1+\chi_s$; then $N_\chi=1$ (as $s\ge 2$), and $(0,\chi)\in S$. As
$\chi(e_s)\ne 0$, there are natural deformations to which the action of $G$
does not extend.

The same argument applies if $Y$ is non-rigid, by replacing $Pic^0(Y)^s$ with
$Pic^0_T(\Y)^s$, where $\Y\to T$ is the restriction of the Kuranishi family of
$Y$ to the locus where $\xi$ stays of type $(1,1)$.
\qed

\begin{rem} {\rm We can find a $Y$ of arbitrary dimension and an ample class
$\xi$ such that deformations of $Y$ for which $\xi$ stays of type $(1,1)$ are
unobstructed; for instance, by taking $Y$ a product of curves of genus at
least two and $\xi$ the canonical class.}
\end{rem}
\begin{thm} Let $d\ge 2$ be an integer. Given any integer $N$, there exists a
point in the moduli space of manifolds of dimension $d$ with ample canonical
class which is contained in at least
$N$ distinct irreducible components.
\end{thm}
\setcounter{equation}{0}
\Pf Without loss of generality, assume that $N\ge d$. Choose arbitrarily
integers $d_1,\ldots,d_N$, each of them $\ge2$ and such that $d_i|d_{i+1}$.
Let
$(Y,L)$ be as in lemma \ref{keylemma}; then for each $k=1,\ldots,N$ there
exists a component of the moduli containing $X$ and having generic
automorphism group isomorphic to $\Z_{d_1}\times\ldots\times\Z_{d_k}$. Hence
$X$ lies in at least $N$ different irreducible components of the moduli.
\qed

In the case of surfaces, this result gives a strong negative answer to the
open problem (ii) on page 485 of \cite{Ca1}.

\begin{thm} Let $G$ be a finite abelian group, and $d\ge 2$ an integer. Then
there exist infinitely many components $M$ of the moduli space of manifolds
of dimension $d$ with ample canonical class such that
$G_M=G$.
\end{thm}
\Pf  Write $G$ as $\Z_{d_1}\times\ldots\times\Z_{d_k}$, with $d_i|d_{i+1}$. If
$k\ge d$, let
$s=k$; if
$k<d$, let $s=d$ and let $d_{k+1}=\ldots=d_s=d_k$.
Choose $(Y,\xi)$ as in
lemma \ref{keylemma}. Applying the lemma to $(Y,i\xi)$ for $i\ge 1$
gives the claimed result.
\qed

In the case of surfaces, another natural question concerns the cardinality of
the automorphism group. Xiao proved in \cite{Xi1} that if $X$ is a minimal
surface of general type,
$\#G\le 52K_X^2+32$ for all abelian subgroups $G$ of $Aut(X)$; it is not
known whether this bound is sharp, but he gives examples to the effect that
any better bound must still be linear in $K_X^2$. It seems natural to ask if
there is a smaller bound if one replaces $Aut(X)$ by $Aut_{\hbox{\rm
gen}}(X)$, the intersection in $Aut(X)$ of $G_M$, for each irreducible
component $M$ containg $X$.  Notice that in Xiao's examples the generic
automorphism group is obviously smaller, so a better bound should be possible.
We prove here that such a bound cannot be less than linear in $K_X^2$.
\begin{prop} There exists a sequence $S_n$ of minimal surfaces of general
type such that
\begin{enumerate}
\item $k_n=K_{S_n}^2$ tends to infinity with $n$;
\item $S_n$ lies on a unique irreducible component, $M_n$;
\item $\#G_{M_n}> 2^{-4}k_n$.
\end{enumerate}
\end{prop}
\Pf Let $n\ge 2$ be an integer. Apply contruction \ref{construction} with
$s=2$, $d_1=d_2=n$, $Y$ a principally polarized abelian surface with
$NS(Y)\isom
\Z$ and
$\xi$ equal to the double of the class of the principal polarization. Choose
$S_n$ to be a cover branched over divisors $D_i$ whose linear equivalence
classes are generic; then all infinitesimal deformations must be Galois, and
the Kuranishi family of $S_n$ is smooth. So $G_{M_n}$ must contain $\Z_n^2$,
hence $\#G_{M_n}\ge n^2$. On the other hand,
$k_n=16(n-1)^2$. Note
that as we only want to bound $G_M$ from below, we don't need to apply
theorem \ref{mainthm}, which would have forced us to choose as class $\xi$ a
higher multiple of the principal polarization.\qed
\begin{rem} {\rm Using the computation of Chern numbers for construction
\ref{construction}, one can determine where the examples constructed so far
lie in the geography of surfaces of general type. For instance by setting all
$d_i$'s equal to $m$ and letting $s$ and $m$ go to infinity, one gets a
sequence of examples where $K^2/c_2$ tends to $2$ from below.}
\end{rem}

\section{Resolution of singularities}
\begin{rem} \label{resnc} Let $\pi:X\to Y$ be a $(G,I)$-cover with $Y$ smooth
and branch locus with normal crossings. Let $Z\to X$ be a resolution of
singularities; then the exceptional locus of $Z$ has uniruled divisorial
components.
\end{rem}
\Pf The question is local on $Y$, so we can assume that $Y$ is affine and
that the line bundles $L_\chi$ and $\O(D_i)$ are trivial.  Let $G'$ be the
abelian group with $\#I$ generators $e_1,\ldots,e_s$, and relations
$m_ie_i=0$ (where $m_i=\#H_i$).  There exists a smooth $G'$-cover $X'$ of $Y$
branched over the
$D_i$ such that the inertia subgroup of $D_i$ is generated by $e_i$, and such
that the map $V\to Y$ factors via $X$. Let $Z'$ be a resolution of
singularities of the fibre product $Z\times_X X'$; we have a commutative
diagram $$
\begin{array}{ccc} Z' & \to & X'\\
\downarrow & & \downarrow\\ Z&\to &X\\
\end{array}
$$ Let $E$ be an irreducible divisorial component of the exceptional locus of
$Z\to X$; its strict transform $E'$ in $Z'$ must be contracted in $X'$ as
$X'\to X$ is finite. As $X'$ is smooth and $Z'\to X'$ is birational, $E'$
must be ruled by \cite{Ab}, therefore $E$ must be uniruled.
\qed

\begin{lem} Let $\Y\to \De$ be a family of smooth manifolds, $\X\to \Y$ an
abelian cover branched on divisors which are all smooth except $D$, of
branching order $n$, which has local equation $f^nh+tg=0$ with $f$, $t$, $h$,
$g$ local coordinates on $Y$ (and $t$ coordinate on $\De$). Then there exists
a morphism $\tilde\Y\to \Y$ such that:
\begin{enumerate}
\item $\tilde\Y\to \Y$ is a composition of blowups with smooth center;
\item the normalization $\tilde \X$ of the induced cover of $\tilde \Y$ is an
abelian cover of $\tilde \Y$ branched over a normal crossing divisor;
\item the exceptional divisors of $\tilde \X\to \X$ have Kodaira dimension
$-\infty$.
\end{enumerate}
\end{lem}
\Pf We will construct $\tilde\Y$ by successive blowups; a local coordinate
and its strict transform after the blowup will be denoted by the same letter.
At each blowing-up step one checks that the normalization of the last
introduced exceptional divisor has Kodaira dimension $-\infty$ (further
blowups change the situation only up to birational maps).

The strategy of the proof is as follows; each blowup introduces a divisor
which is a $\P^r$ bundle (for $r=1,2$), and we prove that the induced cover
of the generic $\P^r$ has Kodaira dimension $-\infty$. We can assume that the
Galois group coincides with the inertia subgroup $H$ of $D$; if this is not
the case, consider the factorization $\X\to \X/H\to Y$, and note that the map
$\X/H\to Y$ is unramified near generic points of $D$, hence after blowing up
the inverse image of the generic $\P^r$ is an unramified cover, which is
therefore a disjoint union of copies of $\P^r$.

We first prove the result on the locus where $h\ne 0$ (this is all one needs
if $\Y$ is a threefold). By changing local coordinates one can assume $h=1$.
Let $n$ be the order of $H$. We distinguish two cases: $n$ even and $n$ odd.
Let $E_1,E_2,\ldots$ be the subsequent exceptional divisors.

\smallskip
\noindent{\sc Case of $n$ even.} Blow up at each step the singular locus
$t=f=g=0$ and look at the $f$ chart. At the first step one obtains $$
z^n=f^2(f^{n-2}+tg)$$ and the total transform of the branch locus $D$ is
$D+2E_1$. The covering restricted to $E_1$ is the composition of a totally
ramified cover of degree $n/2$ and of a double cover ramified over $D\cap
E_1$ which is (on each $\P^2$ in $E_1$) a (possibly reducible) conic. Hence
the cover of $E_1$ is fibered in two-dimensional quadrics (maybe singular).

At the $k$-th step ($1<k\le n/2$) we have $$ z^n=f^{2k}(f^{n-2k}+tg)$$ and
the total transform of $D$ is $$ D+2E_1+\ldots+2kE_k.$$ Again $D$ cuts out a
(possibly reducible) conic on the $\P^2$ fibration of $E_k$; moreover,
$E_k\cap E_i=\emptyset$ if $i<k-1$, and $E_k\cap E_{k-1}$ is (fibrewise) a
line which is not contained in $D$.

If $\xi$ is a generator of the group $H$, the induced cover of $E_k$ is the
composite of a totally ramified cover and of a
cyclic cover of degree $r$, where $r$ is the cardinality of $H/\< \xi^{2k}\>$;
the
cover is ramified on each $\P^2$ on a conic and on a line. The pairs
(inertia group, character) for the branch divisors correspond, via the
bijection defined in \S 2, to $\xi$ for the conic and to $\xi^{-2}$ for the
line.

The canonical bundle of the cover is (fibrewise) the pullback of a multiple
of a line in $\P^2$, the multiple being $$ -3+2\left(\frac{r-1}{r}\right)
+\left(\frac{r/2-1}{r/2}\right)<0
$$ if $r$ is even and $$
-3+2\left(\frac{r-1}{r}\right)+\left(\frac{r-1}{r}\right)<0$$ if $r$ is odd;
in both cases the anticanonical bundle of the cover is ample and the surface
must be of Kodaira dimension $-\infty$.

\smallskip
\noindent{\sc Case of $n$ odd.}  Start by blowing up the singular locus
$t=f=g=0$. At the first step the total transform of $D$ is $D+2E_1$ and the
cover of $E_1$ is totally ramified (as $2$ is prime with $n$), hence the
cover is again $E_1$. If $2k<n$ the same formulas as before hold; we can
repeat the previous argument where $r$ is necessarily odd.

Look now at the $k=(n-1)/2$ case. The total transform of $D$ is $$
D+2E_1+\ldots+(n-1)E_{(n-1)/2}.$$ The strict transform of $D$ is now smooth;
$E_{(n-1)/2}\cap D$ is fibered in singular conics, and we blow up the
singular locus. The center of this blowup does not meet $E_k$ for
$k<(n-1)/2$, and  $D$ and $E_{(n-1)/2}$ have the same tangent space there.
Therefore after blowing one gets an exceptional divisor $E_{(n+1)/2}$
intersecting both $E_{(n-1)/2}$ and $D$ in the same line. The equation (in
the $g$ chart) becomes $$ z^n=f^{n-1}g^n(f+tg).$$
The cover of $E_{(n+1)/2}$ is a
$\P^1$-bundle ramified on a generic $\P^1$ with opposite characters on the
same divisor, hence when normalizing it splits completely. The components of
the total transform of $D$ are smooth, but they meet non-transversally along
the $\P^1$-bundle $f=g=0$.

We now blow up the locus $f=g=0$ and call the exceptional divisor
$F$; the total transform of $D$ is $$
D+2E_1+\ldots+(n-1)E_{(n-1)/2}+nE_{(n+1)/2}+2nF,$$ and $F$ is a $\P^1$-bundle
over a $\P^1$-bundle. The covering of the generic $\P^1$-fibre of $F$ is
ramified of degree $n$ over two points (corresponding to $F\cap D$ and $F\cap
E_{(n+1)/2}$) with opposite characters, hence is again isomorphic to $\P^1$.

\smallskip In both cases the fact that the divisors are smooth and
transversal can be checked at each step out of the center of the next blowup.

We now work in the neighborhood of a point where $h=0$. If $n$ is even, one can
perform the same blowups as in the previous case and check that the same
arguments work. If $n$ is odd, one can perform the first $(n-1)/2$ blowups
as before. After them, the total transform of $D$ has equation
$f^{n-1}(fh+tg)$. In particular (the strict transform of) $D$ is not smooth
any more; we blow up its singular locus, and get a smooth exceptional divisor
$\bar E$. The total transform of $D$ is $$
D+2E_1+\ldots+(n+1)\bar E$$
and is given (in local equations in the $h$ chart) by $$
f^{n-1}h^{n+1}(f+tg).$$
Let $\xi$ be a generator of $H$; the induced cover of $\bar E$ is
cyclic with group $H/\<\xi^{n+1}\>$, hence it is totally ramified and
therefore of Kodaira dimension $-\infty$, being a $\P^2$-bundle.
We are not done because the divisors $D$ and $E_{(n-1)/2}$
are not transversal along $f=g=t=0$; but now we can apply the previous
blowup procedure again.
\qed

\begin{prop}\label{reslemma} Let $\X\to \Y\to \De$ be an abelian cover,
branched over all smooth divisors except one, which has local equation
$f^mh+tg$, where $f,t,g$ are coordinates and $m$ is the order of branching
(where $t$ is the coordinate on $\De$). Then $\X$ and all its transforms via
an $n$-th root base change admit a resolution of singularities such that the
divisorial components of the exceptional divisor all have Kodaira dimension
$-\infty$.
\end{prop}
\Pf The statement without the base change has already been proved; let
$\tilde \X$ be such a resolution. By Hironaka's resolution of singularities
(\cite{Hi}, p.~113, lines 8--4 from the bottom) we can assume that
$\tilde\X_0$ is a normal crossing divisor. Let now $\rho_n:\De\to \De$ be the
map $t\mapsto t^n$. There is a natural birational mapping $\rho_n^*\tilde
\X\to \rho_n^*\X$; moreover $\rho_n^*\tilde \X$ is a cyclic cover of the
manifold $\tilde \X$ ramified over $\tilde\X_0$, which has normal
crossings, hence by remark \ref{resnc} $\rho_n^*\tilde \X$ has a resolution
such that the
divisorial components of the exceptional divisor all have Kodaira dimension
$-\infty$.
\qed

$$
\begin{array}{lll}
\hbox{\rm Dipartimento di Matematica}
&\phantom{aa}
&\hbox{\rm Dipartimento di Matematica e Informatica}
\\
\hbox{\rm Via Sommarive, 7}&&
\hbox{\rm Via Zanon, 6}
\\
\hbox{\rm I 38050 Povo - Italy}&&
\hbox{\rm I 33100 Udine - Italy}
\\
\hbox{\rm fantechi@itnvax.science.unitn.it}&&
\hbox{\rm pardini@ten.dimi.uniud.it}
\end{array}
$$

\end{document}